\documentclass[12pt]{article}
\usepackage{hyperref}
\usepackage{amsfonts,amsmath,amstext,amssymb}
\usepackage{dsfont}
\usepackage{graphicx}

\newtheorem{thm}{Theorem}

\newtheorem{corollary}{Corollary}
\newtheorem{remark}{Remark}

\newenvironment{demo}{\noindent \textit{Proof.} }{\hfill{$\Box$}}

\title{Critical points of the solutions to the $H_R=H_L$ surface equation}

\author{Alma L. Albujer\\
{\tt \small alma.albujer@uco.es}\\
\and 
Magdalena Caballero\\
{\tt \small magdalena.caballero@uco.es}\\
\and
{\small Departamento de Matem\'aticas, Campus de Rabanales,}\\ {\small Universidad de C\'ordoba, 14071 C\'ordoba, Spain}}

\date{}

\begin{document}


\maketitle

\begin{abstract}
Spacelike surfaces with the same mean curvature in $\mathds{R}^3$ and $\mathds{L}^3$ are locally described as the graph of the solutions to the $H_R=H_L$ surface equation, which is an elliptic partial differential equation except at the points at which the gradient vanishes, because the equation degenerates. In this paper we study precisely the critical points of the solutions to such equation. Specifically, we give a necessary geometrical condition for a point to be critical, we obtain a new uniqueness result for the Dirichlet problem related to the $H_R=H_L$ surface equation and we get a Heinz-type bound for the inradius of the domain of any solution to such equation, improving a previous result by the authors. Finally, we also get a bound for the inradius of the domain of any function of class $\mathcal{C}^2$ in terms of the curvature of its level curves.
\vspace{.3cm}


\end{abstract}

\section{Introduction}\label{intro}

A surface in the $3$-dimensional Lorentz-Minkowski space $\mathds{L}^3$ is said to be spacelike if its induced metric is a Riemannian one. Therefore, a spacelike surface can be endowed with two different Riemannian metrics, the one induced by $\mathds{L}^3$ and the metric inherited from the Euclidean space $\mathds{R}^3$. Consequently, we can consider two different mean curvature functions on a spacelike surface related to the previous Riemannian metrics. Let us denote those functions by $H_L$ and $H_R$, respectively.

A surface in $\mathds{R}^3$ is said to be minimal if $H_R$ vanishes identically. Analogously, a maximal surface is a spacelike surface in $\mathds{L}^3$ such that $H_L\equiv 0$. This terminology comes from the fact that those surfaces locally minimize, or maximize respectively, area among all nearby surfaces sharing the same boundary, see~\cite{Lop} and~\cite{Morgan}.

The study of minimal and maximal surfaces has been a topic of wide interest during the last decades. 
In 1983 O. Kobayashi~\cite{Kobayashi} studied, from a local point of view, surfaces in $\mathds{L}^3$ which are simultaneously minimal and maximal. He showed that they are necessarily open pieces of a spacelike plane or of a helicoid in the region where the helicoid is spacelike. 

In 2017, the authors considered in~\cite{AC} the general situation where a spacelike surface in $\mathds{L}^3$ verifies $H_R=H_L$. Those surfaces are locally the graph over a domain $\Omega$ of the Euclidean plane 
of the solutions to a certain partial differential equation called the $H_R=H_L$ surface equation, which is elliptic except at the points at which the gradient vanishes. Among other results, given a solution $u$ to the $H_R=H_L$ surface equation, the authors obtained a bound for the inradius of $\Omega^*$, the subset of $\Omega$ in which the gradient of $u$  does not vanish, see~\cite[Theorem 8]{AC}. Let us recall that the \textit{inradius} of a set in $\mathds{R}^2$ is defined as the supremum of the radius of the closed discs contained in it. Also, in~\cite[Theorem 9]{AC} they proved that if the Dirichlet problem related to the $H_R=H_L$ surface equation has a solution without critical points, then the solution is unique. Notice that in~\cite{AC} the authors do not use the term inradius, but they used instead the term width, which they define as the double of the inradius.

Motivated by the previous results, our goal in this manuscript is to study the critical points of the solutions to the $H_R=H_L$ surface equation. Specifically, given a solution to the $H_R=H_L$ surface equation, we will obtain a necessary condition for a point to be critical, Theorem \ref{critical-points-mean-curvature}. From which we will derive a new uniqueness result for the Dirichlet problem associated to the $H_R=H_L$ equation, Corollary \ref{Dirichlet}, as well as an improved version of \cite[Theorem 8]{AC}, Theorem \ref{width-HRHL}. Inspired by the proof of \cite[Theorem 8]{AC}, we will finally get a bound for the inradius of the domain of any function of class $\mathcal{C}^2$ in terms of the curvature of its level curves, provided a topological condition on the set of critical points of the function is fulfilled, Theorem \ref{width}. 

It is worth pointing out that, recently, several authors have been interested in the study of surfaces (or hypersurfaces) with $H_R=H_L=0$ in more general ambient spacetimes where the problem of considering two different metrics make sense. Specifically, the case where the ambient space is a product has been considered in~\cite{AAS},~\cite{dLR} and~\cite{Kim}, and the case where the ambient is the Heisenberg group has been studied in~\cite{Shin}. The key for the study of simultaneously minimal and maximal surfaces in~\cite{Kim},~\cite{Kobayashi} and~\cite{Shin} is to prove that the level curves of such surfaces are geodesics. However, such reasoning fails at points with horizontal tangent plane, i.e. at critical points, since level curves can be not well defined at those points. 
Furthermore, in~\cite{dLR} the authors consider different techniques and they explicitly impose the assumption of non-existence of critical points. 

The study of critical points of the solutions to elliptic partial differential equations, as well as the study of its level sets, is an issue of special relevance. Both concepts are intimately connected since critical values of $u$ are those at which the level sets can change its topology. 
Let us mention, for instance, a survey by Magnanini~\cite{Mag} (see also references therein), where three different issues about the critical points of the solutions to the Dirichlet problem related to certain elliptic and parabolic partial differential equations are studied. Specifically, the estimation of the local size of the critical set and its location are studied, as well as the dependence of the number of critical points on the boundary values and the geometry of the domain.

\section{Preliminaries}\label{prelim}

Let $\mathds{L}^3$ be the $3$-dimensional Lorentz-Minkowski space, that is, $\mathds{R}^3$ endowed with the metric
\[\langle,\rangle_L=dx^2+dy^2-dz^2,\]
$(x,y,z)$ being the canonical coordinates in $\mathds{R}^3$. The Levi-Civita connections of the Euclidean space $\mathds{R}^3$ and the Lorentz-Minkowski space $\mathds{L}^3$ coincide, so we will just denote them by $\overline{\nabla}$.

A (connected) surface $\Sigma$ in $\mathds{L}^3$ is said to be a spacelike surface if the metric inherited from $\mathds{L}^3$ is a Riemannian one, which is also denoted by $\langle,\rangle_L$. Given a spacelike surface $\Sigma$, there exists a unique future-directed unit normal vector field $N_L$ on $\Sigma$. The mean curvature function of $\Sigma$ with respect to $N_L$ is defined by
\[
H_L=-\frac{1}{2}(k_1^L+k_2^L),
\]
where $k_1^L$ and $k_2^L$ stand for the principal curvatures of $(\Sigma,\langle,\rangle_L)$.

The same topological surface can also be considered as a surface of the Euclidean space $\mathds{R}^3$. In this case, let us denote the induced metric on $\Sigma$ by $\langle,\rangle_R$. It is well-known that $\Sigma$ admits a unique upwards directed unit normal vector field, $N_R$. The mean curvature function of $\Sigma$ with respect to $N_R$ is defined by
\[
H_R=\frac{1}{2}(k_1^R+k_2^R),
\]
where $k_1^R$ and $k_2^R$ stand for the principal curvatures of $(\Sigma,\langle,\rangle_R)$.

A spacelike surface is locally a graph over a domain of the plane $z=0$, which is usually identified with $\mathds{R}^2$, see~\cite[Proposition 12.1.6]{Lop}. Thus, locally $\Sigma=\Sigma_u$, where
\[
\Sigma_u=\{(x,y,u(x,y)):(x,y)\in\Omega\},
\]
for some domain $\Omega\subseteq\mathds{R}^2$ and some smooth function $u\in\mathcal{C}^\infty(\Omega)$.
It is easy to check that $\Sigma_u$ is a spacelike surface if and only if $|Du|<1$, where $D$ and $|\cdot|$ stand for the gradient operator and the norm in the Euclidean plane $\mathds{R}^2$, respectively. Moreover, with a straightforward computation we get the following expressions for the mean curvature functions,
\[
H_L=\dfrac{1}{2}\mathrm{div}\left(\frac{Du}{\sqrt{1-|Du|^2}}\right) \qquad \mathrm{and} \qquad H_R=\dfrac{1}{2}\mathrm{div}\left(\frac{Du}{\sqrt{1+|Du|^2}}\right).
\]
Consequently, any spacelike surface in $\mathds{R}^3$ satisfying $H_R=H_L$ is locally the graph over a certain domain $\Omega\subseteq\mathds{R}^2$ of a solution to the equation
\[
\mathrm{div}\left(\left(\frac{1}{\sqrt{1-|Du|^2}}-\frac{1}{\sqrt{1+|Du|^2}}\right)Du\right)=0,
\]
such that $|Du|<1$. The previous equation is known as the \textit{$H_R=H_L$ surface equation}, which becomes a quasilinear elliptic partial differential equation, everywhere except at those points at which $Du$ vanishes, see~\cite{AC}.

Given a graph $\Sigma_u$ over a domain $\Omega\subseteq\mathds{R}^2$, let $\Sigma_u^\ast$ be the graph of $u$ over the following open set
\[
\Omega^\ast=\{(x,y)\in\Omega:Du(x,y)\neq 0\}.
\]
For any point $p\in\Sigma_u^\ast$, we can consider its corresponding level curve contained in $\mathds{R}^2$. Let $\tilde{k}(\pi(p))$ denote the curvature of such level curve at the point $\pi(p)$, where $\pi$ denotes the natural projection of $\Sigma_u$ onto $\Omega$. Then, in~\cite[Lemma 7]{AC} it has been proved that given any spacelike graph in $\mathds{L}^3$ such that $H_R=H_L$, for any $p\in\Sigma_u^\ast$ it holds
\begin{equation}\label{ineq:kH}
|H_L(p)|\leq \frac{1}{2\sqrt{2}}|\tilde{k}(\pi(p))|.
\end{equation}

\section{On the critical points of the solutions to the $H_R=H_L$ surface equation}
Let $M$ be a differential manifold and let $u:M\longrightarrow \mathds{R}$ be a function of class $\mathcal{C}^r$ with $r\geq 2$. According to its Hessian, the critical points of $u$ can be classified in two types: non-degenerate critical points, those with non-degenerate Hessian, and degenerate critical points, those with degenerate Hessian, see \cite{BG}. 

In the particular case in which $u$ is a smooth function defined over a domain $\Omega\subseteq \mathds{R}^2$, the classification can be stated in terms of the Gaussian curvature of the graph of $u$ either in $\mathds{R}^3$ or in $\mathds{L}^3$: non-degenerate critical points are those with non-vanishing Gaussian curvature, whereas degenerate critical points are those with vanishing Gaussian curvature.

There is a widely known lemma by M. Morse, see \cite[4.2.12]{BG} and \cite[Lemma 2.2]{Mil}, explaining the local behavior of a function of class $\mathcal{C}^r$ with $r\geq 3$ over a differential manifold around a non-degenerate critical point. As a particular case, for smooth functions over a domain $\Omega\subseteq \mathds{R}^2$, this lemma assures the existence of a chart $(U,\phi)$ centered at the critical point $(x_0,y_0)$ such that 
\begin{equation}\label{localbehaviour}
u\circ\phi^{-1}(x,y)=u(x_0,y_0)\pm x^2 \pm y^2.
\end{equation}








Our first result reads as follows. 

\begin{thm}\label{critical-points-mean-curvature}
Let $u$ be a solution to the $H_R=H_L$ surface equation defined on a domain $\Omega\subseteq\mathds{R}^2$ and let $(x_0,y_0)$ be a point in $\Omega$. Then, $H_R(x_0,y_0,u(x_0,y_0))=0$ for any critical point $(x_0,y_0)$ of $u$ in $\Omega$. 
\end{thm}

\begin{demo}
If $(x_0,y_0)$ is a degenerate critical point, then the Gaussian curvature of $\Sigma_u$ at $p$ in $\mathds{R}^3$ vanishes. In \cite[Theorem 4]{AC} the authors proved that given a spacelike surface with $H_R=H_L$, if the Gaussian curvature in $\mathds{R}^3$ vanishes at a point, then the mean curvature also vanishes at that point. And so, $H_R(x_0,y_0,u(x_0,y_0))=0$.

Assume $(x_0,y_0)$ is a non-degenerate critical point. In \cite[Theorem 5]{AC} it is proved that given a compact spacelike surface with (necessarily) non-empty boundary such that $H_R=H_L$, the surface is contained in the convex hull of its boundary. From this result we deduce that if $u$ admits a closed level curve, then its interior can not be contained in $\Omega$. Hence, there are only two possibilities for \eqref{localbehaviour},
\begin{equation}\label{localbehaviour2}
u\circ\phi^{-1}(x,y)=u(x_0,y_0)+x^2 - y^2 \,\, \text{ and }
\,\, u\circ\phi^{-1}(x,y)=u(x_0,y_0)-x^2 + y^2.
\end{equation}

Let us denote by $O$ the codomain of $\phi$. Its domain, $U$, minus the level set through $(x_0,y_0)$ can be divided into the following four domains

$$\begin{array}{cccc}	
U_1&=&\phi^{-1}\left(\{(x,y)\in O\,:\,x^2-y^2>0, x>0\}\right),&\\
U_2&=&\phi^{-1}\left(\{(x,y)\in O\,:\,x^2-y^2>0, x<0\}\right),&\\
U_3&=&\phi^{-1}\left(\{(x,y)\in O\,:\,x^2-y^2<0, y>0\}\right),& \text{ and }\\
U_4&=&\phi^{-1}\left(\{(x,y)\in O\,:\,x^2-y^2<0, y<0\}\right).&
\end{array}$$

We orient the level curves of $u$ in $U-\{(x_0,y_0)\}$ so that its normal vector field points to the direction in which $u$ decreases. 

We consider an open disc centered at $(x_0,y_0)$, $D_{\varepsilon}=D((x_0,y_0), \varepsilon)\subset U$, for some small enough radius $\varepsilon >0$.
Let us demonstrate that there exists a point in $U_1\cap D_{\varepsilon}$ and another one in $U_3\cap D_{\varepsilon}$ such that the curvature of the level curves through both points do not have the same sign. Figure~\ref{fig:thm1} will be useful in the chain of reasonings below. Denote by $\overline{U_1}$ the closure of $U_1$. We take a segment contained in $\overline{U_1}\cap \left(D_{\varepsilon}-\{(x_0,y_0)\}\right)$ intersecting only once each of the following level curves: $$\phi^{-1}\left(\{(x,y)\in O\,:\,x=y, x>0\}\right)\,\, \text{and} \,\, \phi^{-1}\left(\{(x,y)\in O\,:\,x=-y, x>0\}\right).$$ Those curves, the point $(x_0,y_0)$, and the segment define a compact set. We consider the first level curve in $U_1$ that intersects that set and let $(x_1,y_1)$ be a point in the intersection. The point $(x_1,y_1)$ lies on the segment and the level curve is tangent to it. We repeat the same construction for $U_3$, obtaining $(x_3,y_3)$. Taking into account \eqref{localbehaviour2} and the way in which we have oriented the level curves of $u$ in $U$, we can affirm that at one of these two points the normal vector to the level curve points to the interior of its corresponding compact set, while the curve is not locally contained in it. Whereas the normal vector to the level curve at the remaining point points to the exterior of its corresponding compact set. Therefore $(x_1,y_1)$ and $(x_3,y_3)$ are the desired points. 

\begin{figure}[h]
	\centerline{\includegraphics[width=7cm]{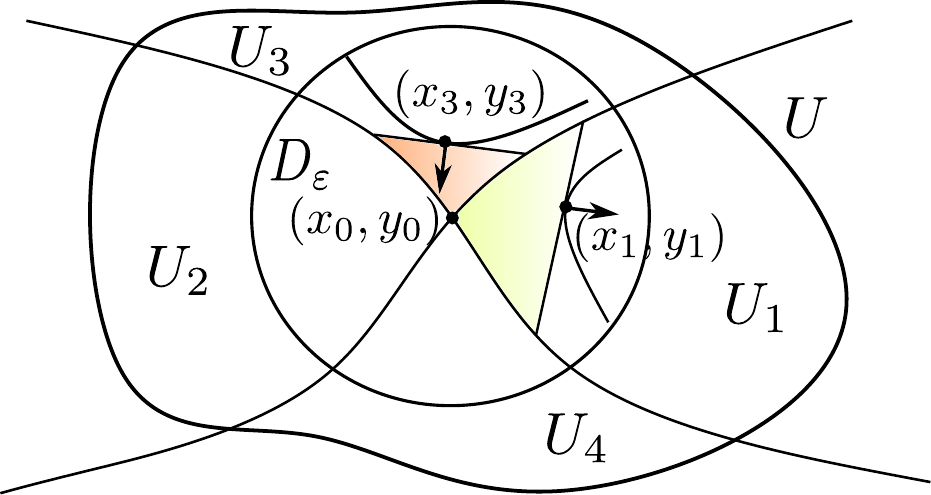}}  \caption{Construction of $(x_1,y_1)$ and $(x_3,y_3)$} 
	\label{fig:thm1}
\end{figure}  

The curvature of the level curves is a continuous function in $D_{\varepsilon}-\{(x_0,y_0)\}$. As such, there exists a point $(\tilde{x},\tilde{y})$ in the punctured disc whose level curve has vanishing curvature at that point. From inequality~\eqref{ineq:kH} $H_L(\tilde{x},\tilde{y},u(\tilde{x},\tilde{y}))=0$. Taking limit when $\varepsilon$ approaches zero, we complete the proof. 

\end{demo}


We have mentioned in the Introduction that if the Dirichlet problem related to the $H_R=H_L$ surface equation has a solution without critical points, then the solution is unique, see \cite[Theorem 9]{AC}. Therefore as an immediate corollary of Theorem~\ref{critical-points-mean-curvature} we get a new uniqueness result for the Dirichlet problem associated to the $H_R=H_L$ surface equation. 

\begin{corollary}\label{Dirichlet}
If $u$ is a solution to the Dirichlet problem associated to the $H_R=H_L$ surface equation over a domain $\Omega\subseteq\mathds{R}^2$ such that its graph has non-vanishing mean curvature, then $u$ is the only solution.
\end{corollary}

In 1955 E.~Heinz used the classical divergence theorem to prove that given a graph in $\mathds{R}^3$ defined over a disk of radius $R$ in $\mathds{R}^2$ centered at the origin, if $|H_R|\geq c>0$ for a certain constant $c$, then $R\leq \dfrac{1}{c}$, see \cite{Heinz}. 
This inequality can be restated in terms of the inradius as follows, let $u$ be a smooth function over a domain $\Omega\subseteq \mathds{R}^2$, then $$\mathrm{inradius}(\Omega)\leq \dfrac{1}{\mathrm{inf}|H_R|}.$$  

In \cite[Theorem 8]{AC} the authors got the following bound for the graphs satisfying $H_R=H_L$, $$\mathrm{inradius}(\Omega^*)\leq \dfrac{1}{2\sqrt{2}\,\mathrm{inf}|H_R|},$$ where, as we have mentioned before, $\Omega^*$ is the set of points at which $Du$ does not vanish.  Combining this result with Theorem 1 we get.

\begin{thm}\label{width-HRHL}
Let $\Sigma_u$ be a spacelike graph in $\mathds{L}^3$ over a domain $\Omega\subseteq\mathds{R}^2$ such that $H_R=H_L$. Then 
$$\mathrm{inradius}(\Omega)\leq \dfrac{1}{2\sqrt{2}\,\mathrm{inf}|H_R|}.$$
\end{thm}







The following corollary is an immediate consequence of the previous result.


\begin{corollary}
Let $\Sigma_u$ be a spacelike graph in $\mathds{L}^3$ defined over a domain of infinite inradius such that $H_R=H_L$. Then $\mathrm{inf}|H_R|=0$. 

Equivalently, there do not exist spacelike graphs over an infinite inradius domain satisfying $H_R=H_L$ and $|H_R|\geq c$ for certain positive constant $c$.
\end{corollary}

\section{A bound for the inradius of a graph depending on its level curves}

Inspired by the proof of~\cite[Theorem 8]{AC}, we get our next theorem. 

\begin{thm}\label{width}
Let $u$ be a $\mathcal{C}^2$ function defined over a domain $\Omega\subseteq\mathds{R}^2$. Let $\tilde{k}$ denote the curvature of its level curves. Denote by $A$ the set of critical points of $u$ and by $A'$ its accumulation set. Define $A^{\left.2\right)}=(A')'$, $A^{\left. n\right)}=(A^{\left.n-1\right)})'$ for $n>2$. If there exists $n$ such that $A^{\left. n\right)}=\emptyset$, then $$\mathrm{inradius}(\Omega)\leq \dfrac{1}{\mathrm{inf}|\tilde{k}|}.$$
\end{thm}





\begin{demo}
Assume $\mathrm{inf}|\tilde{k}|\neq 0$. We consider all the level curves in $\Omega$, and we orient them so that the normal vector field points to the direction in which $u$ decreases. 

We begin with the case $A'=\emptyset$. 

Let us assume $\mathrm{inradius}(\Omega)>1/\mathrm{inf}|\tilde{k}|$. Hence, for each $c$ such that 
\begin{equation}\label{ineq-kc}
1/\mathrm{inradius}(\Omega)<c<\mathrm{inf}|\tilde{k}|,
\end{equation} there exists a closed disc with center at a point $q\in\Omega$ and radius $1/c$, $\bar{B}_q(1/c)$, contained in $\Omega$.

Since $A'=\emptyset$, there are only a finite number of critical points in the disc. If necessary, we take a bigger $c$ so that there is no critical point on the boundary of the disc. 

For each critical point in $B_q(1/c)$ we take a small enough open disc centered at it such that their closures do not intersect and are contained in $B_q(1/c)$. We denote by $D$ the compact set obtained by subtracting those discs to $\bar{B}_q(1/c)$. 

Our function $u$ has no critical points in $D$, therefore it attains its extremal values on its boundary. Let us assume first that the maximum $p$ is attained at a point on the boundary of $\bar{B}_q(1/c)$. In that case, the level curve through $p$ lies in $\Omega\setminus B_q\left(1/c\right)$. And so, it is tangent to the boundary of the disc at $p$. The normal vector to the curve at $p$ points to the interior of the disc, while the curve is not locally contained in it. Consequently, inequality \eqref{ineq-kc} implies that $\tilde{k}< -c$ at $p$. Otherwise, the maximum is attained at a point on the boundary of one of the open discs centered at a critical point. In this case, the normal of the level curve through that point is directed to the exterior of the disc, whereas the level curve is contained in its closure. We conclude that, no matter where the maximum value is attained, the curvature of the level curve at that point is always negative, see Figure~\ref{fig:thm3}. In an analogous way we prove that given a point at which $u$ attains a minimum, the curvature of the level curve at that point is positive. By a continuity argument we conclude that there exists a point in $D$ at which the curvature of the level curve through the point vanishes, which is a contradiction.

\begin{figure}[h]
	\centerline{\includegraphics[width=13cm]{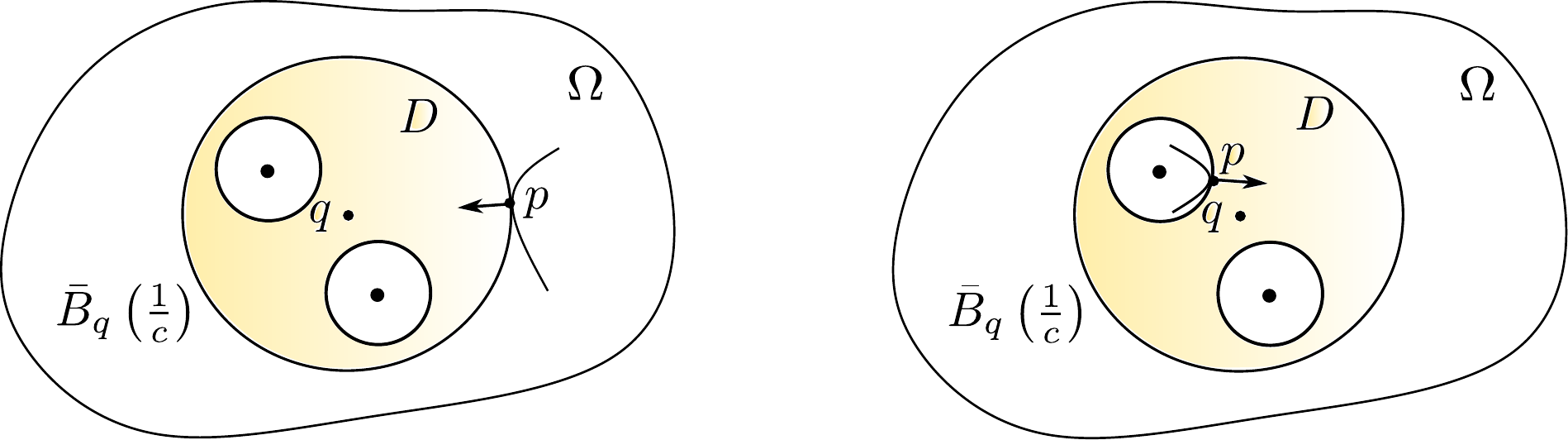}}  \caption{Level curve at the maximum $p$} 
	\label{fig:thm3}
\end{figure}  

To complete the proof we only need to prove that if there exists $n>1$ such that $A^{\left. n\right)}=\emptyset$, we can always choose $c$ and a finite number of open discs such that their closures do not intersect and are contained in $B_q(1/c)$, and satisfying that no critical point lies on the set obtained by subtracting those discs to $\bar{B}_q(1/c)$, which will be called $D$.


If $A^{\left. n\right)}=\emptyset$ and $A^{\left. n-1\right)}\neq\emptyset$, there is only a finite number of points of $A^{\left. n-1\right)}$ in $\bar{B}_q(1/c)$, otherwise they will accumulate. Choose $c$ big enough such that none of them lie on the boundary of $\bar{B}_q(1/c)$. Take an open disc centered at each of those points such that their closures do not intersect and are contained in $B_q(1/c)$. We denote by $D_{n-1}$ the set obtained by subtracting those discs to $\bar{B}_q(1/c)$. In $D_{n-1}$ there is only a finite number of points of $A^{\left. n-2\right)}$, otherwise they will accumulate. Take $c$ and the previously chosen open discs big enough such that neither of the points in $A^{\left. n-2\right)}$ lie on the boundary of $D_{n-1}$. Take an open disc centered at each of those points such that their closures do not intersect and are contained in the interior of $D_{n-1}$. We denote by $D_{n-2}$ the set obtained by subtracting those discs to $D_{n-1}$. In a finite number of steps we construct $D=D_0$.

\end{demo}

\begin{remark}
In the previous result we can substitute the set of the critical points of $u$ for the set of points at which the level set is not locally a curve.
\end{remark}

\bibliographystyle{amsplain}

\end{document}